\documentclass[a4paper,12pt]{article}

\usepackage[T2A]{fontenc}
\usepackage[utf8]{inputenc}   
\usepackage[russian]{babel}

\usepackage{CJKutf8}

\usepackage{amsmath,amstext,amsfonts,amsthm,amssymb}
\usepackage{eucal}
\usepackage{graphicx}
\renewcommand{\subsubsection}[1]{\addtocounter{subsubsection}{1}
{\ \\[3pt]\bf \thesubsubsection. \  #1} }

\swapnumbers

{  \theoremstyle{definition}

}
\newcommand{\Coker}{\operatorname{Coker}}
\newcommand{\Cone}{\operatorname{Cone}}

\newcommand{\Ker}{\operatorname{Ker}}
\newcommand{\Loc}{\operatorname{Loc}}

\newcommand{\Perv}{\operatorname{Perv}}

\newcommand{\hra}{\hookrightarrow}
\newcommand{\iso}{\overset{\sim}{\longrightarrow}}
\newcommand{\isom}{\overset{\sim}{=}}
\newcommand{\lra}{\longrightarrow}

\newcommand{\lla}{\longleftarrow}

\newcommand{\bea}{\begin{eqnarray*}}
\newcommand{\eea}{\end{eqnarray*}}
\newcommand{\bean}{\begin{eqnarray}}
\newcommand{\eean}{\end{eqnarray}}

\newcommand{\CA}{\mathcal{A}}
\newcommand{\CB}{\mathcal{B}}
\newcommand{\CC}{\mathcal{C}}

\newcommand{\CM}{\mathcal{M}}

\newcommand{\CW}{\mathcal{W}}

\newcommand{\BC}{\mathbb{C}}

\newcommand{\BF}{\mathbb{F}}
\newcommand{\BG}{\mathbb{G}}

\newcommand{\BZ}{\mathbb{Z}}

\newcommand{\nc}{\newcommand}

\nc{\Id}{\text{Id}}
\nc{\la}{\lambda}

\begin{document}

\centerline{\bf WITT RING LOCALIZED}

\


\bigskip\bigskip

\centerline{Vadim Schechtman}

\bigskip\bigskip

\begin{CJK}{UTF8}{min}


\end{CJK}


 \centerline{October 29, 2020}

\

\


\centerline{\bf Introduction}

\

The rings of $p$-typical Witt vectors are interpreted as spaces of vanishing  cycles 
for some perverse sheaves over a disc. 

This allows to "localize"\ an isomorphism emerging in Drinfeld's theory of prismatization [Dr], Prop. 3.5.1, namely to  express it  as "an integral"\  of a standard  
exact triangle on the disc.


I am grateful to A.Beilinson, V.Drinfeld, J.Tapia and B.Toen for  consultations.

\

\

\centerline{\bf \S 1. Some elementary examples of perverse sheaves over $(\BC, 0)$}

\

{\bf 1.1. Perverse sheaves over the disc.} Let $\CA$ be an abelian category. 
We denote by $\Perv(D, 0; \CA) = \Perv(\BC, 0; \CA)$ the category 
of diagrams $\CM$ in $\CA$
$$
\Phi(\CM)\begin{matrix} v\\ \lra \\ \lla \\ u
\end{matrix} \Psi(\CM)
\eqno{(1.1.1)}
$$
where $u, v$ are such that

$(Inv)$ $T_\Psi:= 1 - vu$ is invertible. 

The last condition is equivalent to:

$(Inv')$ $T_\Phi:= 1 - uv$ is invertible.

We use the notations
$\Perv(\BC, 0) = \Perv(\BC, 0; \CA b)$, etc.,  
where $\CA b$ is the category of abelian groups.

\

(a) {\it Cohomology}

\

By definition
$$
\Gamma_c(D;\CM): \ 0 \lra \Phi(\CM) \overset{v}\lra \Psi(\CM) \lra 0,
$$
a complex living in degrees $0, 1$, and 
$$
\Gamma(D;\CM): \ 0 \lra \Psi(\CM) \overset{u}\lra \Phi(\CM) \lra 0,
$$
a complex living in degrees $-1, 0$.

\

(b) {\it Standard functors}

\

Let $\Loc(\BC^*; \CA)$ denote the category whose objects are diagrams in $\CA$
$$
\Psi\overset{T}\iso \Psi
$$
We have a  forgetful functor 
$$
j^*: \ \Perv(\BC, 0; \CA) \lra \Loc(\BC^*; \CA),\ (\Phi, \Psi, u, v)\mapsto (\Psi, 1 - vu)
$$
which admits a left adjoint
$$
j_!(\Psi, T) = (\Psi, \Psi, 1, 1 - T)
$$
and a right adjoint 
$$
j_*(\Psi, T) = (\Psi, \Psi, 1 - T, 1)
$$

\

(c) {\it Exact triangles}

\

We have the usual exact triangles
$$
i_!i^!\CM \lra \CM\lra j_*j^*\CM,
$$
so that
$$
i^!\CM = \biggl(\Phi \overset{v}\lra \Psi\biggr) = \Gamma_c(D;\CM),
$$
in degrees $0, 1$, and
$$
j_!j^!\CM \lra \CM \lra i_*i^*\CM,
$$
so that
$$
i^*\CM = \biggl(\Psi \overset{u}\lra \Phi\biggr) = \Gamma(D;\CM)
$$
in degrees $-1, 0$, as it should be.

\

{\bf 1.2. Example.} Let $A$ be a ring, $b\in A$. For $n, m\in \BZ_{\geq 0}$ define 
an object 
$$
M(A, b)_{m, n} \in \Perv(\BC, 0)
$$ 
by a diagram
$$
A/(b^{n})\begin{matrix} v_{m,n}\\ \lra \\ \lla \\ u_{m,n}
\end{matrix} A/(b^{n+m}) 
\eqno{(1.2.1)}
$$
where $v = v_{m,n}$ is the multiplication by $b^m$ and $v = u_{m,n}$ is the canonical projection.

The composition $v_{m,n}u_{m,n} = b^m$ is nilpotent on $A/(b^{n+m})$, therefore
$$
T_{m,n} = 1 - v_{m,n}u_{m,n}
$$
is invertible.

We have obvious maps in $\Perv(\BC, 0)$
$$
M(A, b)_{m, n} \lra M(A, b)_{m - 1, n};
$$
passing to the inverse limit $\lim_{\overset{n}\leftarrow}$ we get a perverse sheaf
$$
\hat M(A, b)_m: \hat A_b \begin{matrix} b^m\\ \lra \\ \lla \\ 1
\end{matrix} \hat A_b
$$
which is the $j_!$ of a local system on $\BC^*$
$$
\hat M(A, b)_m = j_!(\hat A_b, 1 - b^m)
$$
where $j:\ \BC^*\hra \BC$ (cf. 1.1 (b)).

Similarly, a left $A$-module  $N$ gives rise to  an object 
$$
M_{m,n}(N, b) = (N/b^{n}N, N/b^{n+m}N; b^m, 1)\in \Perv(\BC, 0) 
$$

{\bf 1.3. Example.} Take $A = \BZ, \ b = p$ a prime number, $m = 1$. We get a perverse sheaf
$$
\hat M(\BZ, p)_1: \hat \BZ_p \begin{matrix} p\\ \lra \\ \lla \\ 1
\end{matrix} \hat \BZ_p
$$
Note that $\hat \BZ_p = W(\BF_p)$ (the ring of $p$-typical Witt vectors), with $F = 1$ and $V = p$. 

We will see generalizations of this example in the next Section.

{\bf 1.4. Ring objects and modules.} Let $\CA$ be an abelian tensor category. 

Let us call an object
$B =(B_0, B_1, u, v)\in \Perv(\BC, 0; \CA)$ {\it a ring object} if

(i) $B_0, B_1$ are rings in $\CA$;

(ii) $u$ is a ring homomorphism, i.e.
$$
u(xy) = u(x)u(y)
$$
(iii) $v$ satisfies  a projection formula
$$
v(u(y)x) = y v(x).
$$

An object $N = ((N_0, N_1, u, v)\in \Perv(\BC, 0; \CA)$ is a (left) {\it $B$-module} 
if  

(iv) $N_i$ is a left $B_i$-module, $i = 0, 1$;

(v) 
$$
u(a x) = u(a)u(x),\ a\in B_1, x\in M_1;
$$

(vi)
$$
v(b y) = v(b)v(y),\ b\in B_0, y\in M_0.
$$

In the examples 1.2 above, $M_{m,n}(A, b)$ are ring objects, and 
$N_{m,n}(M, b)$ are left $M_{m,n}(A, b)$-modules.

\

\centerline{\bf \S 2. Examples related to the Witt ring}

\

Below there are some more examples.

{\bf 2.1. A perverse sheaf $\CW$ from the Witt ring: "variation = Verschiebung".} Let $p$ be a  
prime number. 

Let $A$ be a commutative $\BF_p$-algebra, and  $W = W_p(A)$ its 
ring of $p$-typical Witt vectors. 

Recall that $W$ comes equipped with two self-maps:
$$
F:\ W\lra W,
$$
"Frobenius"\ , and 
$$
V:\ W\lra W,
$$
"Verschiebung".

$F$ is a ring homomorphism, whereas $V$ satisfies the projection formula 
$$
V(F(x)y) = xV(y)
$$
Also
$$
V\circ F = F\circ V = p
$$
Consider the diagram
$$
\CW:\ W \begin{matrix} V\\ \lra \\ \lla \\ F
\end{matrix} W
\eqno{(2.1.1)}
$$
(so, $v = V$, $u = F$). 

Note that
$$
T: = 1 - vu = 1 - FV = 1 - p,
$$
or $p\in \BZ_p\subset W$ is pro-nilpotent, therefore $T$ is invertible. 

It follows that $\CW\in \Perv(\BC, 0)$; it is a ring object. 

We can interchange $F$ and $V$ and define the "Fourier dual" sheaf $\CW^\vee$.

For $A = \BF_p$ we get the sheaf from 1.3.

\

{\bf 2.2. From Dieudonn\'e modules.} More generally, let $M$ be a finitely generated  topological 
$W$-module (for example a $W$-module of finite length) equipped with two operators $F, V:\ M\lra M$ such that
$$
VF = p\cdot \Id_M,
$$
cf. [D], p. 69.

These data define a perverse sheaf
$$
\CM: M \begin{matrix} V\\ \lra \\ \lla \\ F
\end{matrix} M
\eqno{(2.2.1)}
$$
Indeed, since $M$ is finitely generated, the operator $p$ is pro-nilpotent on $M$, 
whence $T:= 1 - VF$ is invertible.  

For example, to a $p$-divisible group $G$ there corresponds a perverse sheaf 
$\CM(G)$ associated to its Dieudonn\'e module $M(G)$.

\


\centerline{\bf \S 3. Some objects related to prismatization.}

\

{\bf 3.1.} Let $k$ be a commutative algebra, and $W(k) = W_p(k)$. Note that $FV = p$.
We suppose that $W(k)$ is $p$-adically complete (for example a field of char $p$, or 
$\BZ/p^n\BZ$, etc.)
 

Consider a perverse sheaf
$$
\CW^\vee(k):\ W(k) \begin{matrix} F\\ \lra \\ \lla \\ V
\end{matrix} W(k)
\eqno{(3.1.1)}
$$
Following [Dr], [B] consider
$$
\Ker(W(k)\overset{F}\lra W(k))
$$
In our interpretation it is nothing else than $H^{-1}(D;\CW^\vee(k))$. 

Lemma 3.2.6 from [D] says that this group is isomorphic to $\BG_a^\sharp(k)$.


\

{\bf 3.2.} Let us apply to $\CW^\vee(k)$ the formalism of standard functors from 1.1.

We have
$$
j_!j^*\CW^\vee(k):\ W(k) \begin{matrix} p\\ \lra \\ \lla \\ 1
\end{matrix} W(k)
$$
and the adjointness map looks as follows:
$$
j_!j^*\CW^\vee(k) \lra \CW^\vee(k) = (V, 1).
$$
Its cone is a complex concentrated at the origin $0\in D$:
$$
i_*i^*\CW^\vee(k):\ \Cone(V) \begin{matrix} 0\\ \lra \\ \lla \\ 0
\end{matrix} 0
$$
Note that $V: W(k)\lra W(k)$ is injective and its cokernel is $k = \BG_a(k)$, so 
$$
\Cone(V: W(k)\lra W(k)) = \BG_a(k).
$$
Let us denote
$$
\CC_F(k) := \Cone(F:\ W(k)\lra W(k)),\ 
$$
$$
\CC_p(k) := \Cone(p:\ W(k)\lra W(k))
$$
We have $FV = p$, whence a map
$$
(V, 1):\ \CC_p(k)\lra \CC_F(k)
$$
If we apply the functor $\Gamma_c(D, ?)$ to the exact triangle 
$$
\CW^\vee(k) \lra i_*i^*\CW^\vee(k) \lra j_!j^*\CW^\vee(k)[1]
\eqno{(3.2.1)}
$$
we get an exact triangle
$$
\CC_F(k)[-1] \lra \BG_a(k) \lra \CC_p(k).
\eqno{(3.2.2)}
$$
The corresponding cohomology exact sequence is 
$$
0 \lra \Ker(p) \lra \Ker F \lra \BG_a(k) \lra \Coker(p) \lra \Coker F \lra 0
$$

\

Let us suppose that $F$ is surjective, for example this is true if $k$ is a perfect field of 
char $p$. In that case
$$
\CC_F(k) \isom \Ker F[1] \isom \BG_a^\sharp(k)[1],  
$$
and (3.2.2) can be rewritten as an isomorphism
$$
\Cone(\BG_a^\sharp(k) \lra \BG_a(k)) \isom \CC_p(k) 
$$
which is nothing else but the isomorphism from [Dr], Proposition 3.5.1, on the level of 
$k$-points.

{\bf 3.3.} Otherwise, we can replace our target category $\CA$ by a category $\CB$ of fpqc sheaves 
on the category of $p$-adic schemes, and consider $W, \Ker F$, etc. as 
objects of $\CB$. 

In that case 
$$
F:\ W \lra W
$$
is surjective\footnote{I am grateful to B.Toen for this remark}, so
$$ 
\Gamma_c(D, \CW^\vee) \isom \Ker F[-1]\isom \BG_a^\sharp[-1]
$$
The standard exact triangle (3.2.1) takes a form
$$
(E):\ \CW^\vee\lra i_*i^*\CW^\vee\lra j_!j^*\CW^\vee[1],
\eqno{(3.3.1)}
$$
and its "integral"\ $\Gamma_c(D;E)$  will be isomorphic to an exact triangle
$$
\Gamma_c(D;E):\ \BG_a^\sharp \lra \BG_a \lra \Cone(p: \CW^\vee \lra \CW^\vee)
$$
which coincides with the assertion [Dr], Proposition 3.5.1.


\bigskip\bigskip

\centerline{\bf References}

\bigskip\bigskip

[B] A.Beilinson, Drinfeld's approach to the prismatic cohomology, a lecture in MPhTI, 2019, see on youtube


[D] M.Demazure, Lectures on $p$-divisible groups

[Dr] V.Drinfeld, Prismatization, arxiv:2005.04746

\newpage



\end{document}